\documentclass[12pt]{amsart}

\usepackage{amsmath, amssymb}
\usepackage{mathtools}
\usepackage{bbm}
\usepackage[utf8]{inputenc}
\usepackage{calrsfs}
\DeclareMathAlphabet{\pazocal}{OMS}{zplm}{m}{n}
\usepackage{color}
\usepackage{overpic}
\usepackage{hyperref}
\usepackage{soul}
\usepackage{dsfont}
\usepackage{enumitem}

\usepackage[english]{babel}

\usepackage{array}

\usepackage{caption}
\makeatletter
\@namedef{subjclassname@2020}{\textup{2020} Mathematics Subject Classification}
\makeatother

\usepackage{float}

\newtheorem{theorem}{Theorem}[section]
\newtheorem{lemma}[theorem]{Lemma}
\newtheorem{prop}[theorem]{Proposition}

\theoremstyle{definition}

\newtheorem{propdefi}[theorem]{Proposition and Definition}
\newtheorem{cordefi}[theorem]{Corollary and Definition}

\theoremstyle{remark}
\newtheorem{remark}[theorem]{Remark}

\numberwithin{equation}{section}

\allowdisplaybreaks

\begin{document}

\title[Discrete Hyperbolic Cosecant Distributions]{Discrete Hyperbolic Cosecant Distributions}


\author{Leonard J. E. Pleschberger}
\address{Leonard Jobst Eberhard Pleschberger, Neptunstra\ss e 18, D\"usseldorf D-40223 \& Sechshauser Stra\ss e 24/1/15, A-1150 Vienna.}
\email{leonard.pleschberger@gmail.com}

\date{\today}

\begin{abstract}
We introduce three families of discrete hyperbolic cosecant probability distributions, supported either bilaterally on $\mathbb{Z}$ or unilaterally on $\mathbb{N}_0$, inspired by series evaluations of Srinivasa Ramanujan. We determine their exact normalization constants in terms of singular values of complete elliptic integrals and the Gau\ss ian constant $G=\varpi/\pi$ with the lemniscate constant $\Gamma^2(1/4)/(2\sqrt{2\pi})$ by explicitly calculating some series, e.\,g.
$$
\sum_{k \in \mathbb{Z}} 2k \operatorname{cosech}(k \pi) = G^2 + \frac{1}{\pi}.
$$
Furthermore, by connecting these distributions to Fourier transforms of powers of the hyperbolic secant and applying the Poisson summation formula, we obtain closed-form evaluations for even moments and hyperbolic secant lattice sums, e.\,g.
$$
\sum_{k \in \mathbb{Z}} \mathrm{sech}^4(k \pi) = \frac{2}{3} \left(G^4 + G^2 + \frac{1}{\pi}\right)
$$
as polynomials of the Gau\ss ian constant.
\end{abstract}

\maketitle

\section{Introduction}

We introduce three types of discrete hyperbolic cosecant probability distributions, denoted either by $\mathsf{Cosech}$ in the bilateral case supported on $\mathbb{Z}$, or by $\underline{\mathsf{Cosech}}$ in the unilateral case supported on the non-negative natural numbers $\mathbb{N}_0$; this family of distributions is inspired by some series identities provided by Srinivasa Ramanujan. These distributions offer a view of the classical results through the lenses of probability theory and lemniscatic geometry.\\

Ramanujan's formulae in his second notebook have a very general setup. They contain terms of Jacobi's complete elliptic integral of the first and second kind, $K$ and $E$, viz.
$$
K(k) := \int_0^{\pi/2} \frac{\mathrm{d}\vartheta}{\sqrt{1-k^2 \sin^2(\vartheta)}}, \quad E(k) := \int_0^{\pi/2} \sqrt{1-k^2 \sin^2(\vartheta)} \,\mathrm{d}\vartheta.
$$
Here, $k\in(0,1)$ is called the elliptic modulus with complementary modulus $k'$, such that $k^2+k'^2=1$, cf. \cite[p.\,31]{jacobi_1829}. For $r>0$, let $k_r \in (0,1)$ be the unique elliptic modulus with $k'_r:=\sqrt{1-k^2_r}$, $K_r := K(k_r)$, and $K'_r = K(k'_r)$, such that $K'
_r/K_r = \sqrt{r}$ holds. The elliptic lambda-star function then yields
$\lambda^*(r) = k_r$. According to \cite{selberg_1967}, in case of $r \in \mathbb{Q}^+$, the corresponding elliptic integrals $K_r$ and $K'_r$ can be expressed in closed form using a finite number of values of the gamma function. For $r \in \mathbb{N}$, $k_r$ is called an elliptic singular modulus and the corresponding elliptic integrals $K_r$ can often be calculated explicitly. Cf. \cite[p.\,298]{borwein_1987} for the values of $K_r$ for the singular elliptic moduli $k_1, \dots, k_{16}$; numerous further values of singular elliptic moduli are stated in the exercises in the same book. For the calculation of $E_r := E(k_r)$, we use the elliptic alpha function $\alpha(r)$ which is implicitly defined via Equation (5.1.2) in \cite[p.\,152]{borwein_1987}; it yields
\begin{equation}\label{eq:E_r_alpha}
E_r = \frac{\pi}{4K_r\sqrt{r}} + \left(1-\frac{\alpha(r)}{\sqrt{r}}\right)K_r
\end{equation}
and the Tables 5.2a., 5.2b., and 5.3. on p. 172--173 in the same book provide many special values for $\alpha(r)$ for $r \in \mathbb{Q}^+$. For elliptic singular values $k_r$, Ramanujan's notation in the Chapters 17 and 18 of his second notebook can be translated into the common language via $x=k^2_r$, $y=\sqrt{r} \pi$, and $z=2 K_r / \pi$, cf. Entry 6 of Chapter 17 in \cite[p.\,102]{ramanujan_1991}.

Let $\varpi := \Gamma^2(1/4)/(4\sqrt{\pi})$ denote the lemniscate constant. A particularly important elliptic singular modulus occurs in the lemniscatic case where $r=1$. Then, we get $K'_1 = K_1 = \varpi/\sqrt{2}$ and $\lambda^*(1)= k_1 = k'_1 = 1/\sqrt{2}$. Hence, in the lemniscatic case, we may substitute $x=1/2$, $y=\pi$ and $z=\sqrt{2}G$ into Ramanujan's equations. Finally, Gau\ss’ constant $G:=\varpi/\pi$ yields particularly simple expressions for $K_r$, e.\,g. for its values at singular elliptic moduli $k_r$ for $r=s^2$ with $s \in \{1,\dots,5\}$. In this case, also the normalization constants for our discrete type I and type III hyperbolic cosecant distribution are given in notably simple closed forms.\\

\renewcommand{\arraystretch}{2}
\begin{table}[!tbh]
    \centering
    \begin{tabular}{|c|c|c|c|}
\hline
Type & Distribution & Normalizing Constant & Weight Function\\
\hline
I & $\mathsf{Cosech} = \mathsf{Cosech}_{X_1}$ & $\frac{\pi}{\pi G^2 + 1}$ & $ 2\kappa \operatorname{cosech}\left(\kappa \pi \right)$\\
I & $\mathsf{Cosech}_{Y_1}$ & $ \frac{1}{\sqrt{2} G^2}$ & $(2\kappa+1) \operatorname{cosech} \left((2\kappa+1) \frac{\pi}{2}\right)$\\
I & $\mathsf{Cosech}_{Z_1}$ & $\frac{\pi}{\pi G^2+1 + \sqrt{2}\pi G^2} $ & $ \kappa \operatorname{cosech} \left(\kappa \frac{\pi}{2} \right)$\\
\hline
I & $\mathsf{Cosech}_{X_r}$ & $\frac{\sqrt{r} \pi^2}{\pi + 4 \alpha(r) K^2_r}$ & $ 2\kappa \operatorname{cosech}\left(\sqrt{r} \kappa \pi \right)$\\
I & $\mathsf{Cosech}_{Y_r}$ & $\frac{\pi^2}{4 k_r K^2_r}$ & $(2\kappa+1) \operatorname{cosech} \left(\sqrt{r} (2\kappa+1) \frac{\pi}{2}\right)$\\
I & $\mathsf{Cosech}_{Z_r}$ & $\frac{\sqrt{r} \pi^2}{\pi + 4\alpha(r)K^2_r + 4 \sqrt{r} k_r K^2_r}$ & $ \kappa \operatorname{cosech} \left(\sqrt{r} \kappa \frac{\pi}{2} \right)$\\
\hline
I & $\underline{\mathsf{Cosech}}_{X_1}$
& $\frac{2\pi}{\pi G^2 + 3}$
& $(2n + 2) \operatorname{cosech}\left((2n + 2) \frac{\pi}{2} \right)$\\
I & $\underline{\mathsf{Cosech}}_{Y_1}$
& $\frac{\sqrt{2}}{G^2}$
& $(2n + 1) \operatorname{cosech}
\left((2n + 1) \frac{\pi}{2}\right)$\\
I & $\underline{\mathsf{Cosech}}_{Z_1}$
& $\frac{2\pi}
{\left(1 + \sqrt{2}\right) \pi G^2 + 3}$
& $(n + 1) \operatorname{cosech}
\left((n + 1) \frac{\pi}{2} \right)$\\
\hline
I & $\underline{\mathsf{Cosech}}_{X_r}$
& $\frac{2 \sqrt{r} \pi^2}{4 \alpha(r) K^2_r + 3\pi}$
& $(2n+2) \operatorname{cosech}\left(\sqrt{r} (2n + 2) \pi \right)$\\
I & $\underline{\mathsf{Cosech}}_{Y_r}$
& $\frac{\pi^2}{2 k_r K^2_r}$
& $(2n + 1) \operatorname{cosech}
\left(\sqrt{r} (2n + 1) \frac{\pi}{2}\right)$\\
I & $\underline{\mathsf{Cosech}}_{Z_r}$
& $\frac{2\sqrt{r} \pi^2}
{4\alpha(r)K^2_r + 4 \sqrt{r} k_r K^2_r + 3\pi}$
& $(n+1) \operatorname{cosech}
\left(\sqrt{r} (n+1) \frac{\pi}{2} \right)$\\
\hline
\end{tabular}
\caption{Discrete bilateral (lemniscatic) type I $\mathsf{Cosech}$ distributions of the first ($X$), second ($Y$), and third ($Z$) kind for $\kappa \in \mathbb{Z}$ and $r>0$ ($r=1$); further the discrete unilateral (lemniscatic) type I $\underline{\mathsf{Cosech}}$ distributions of the first ($X$), second ($Y$), and third ($Z$) kind for $n \in \mathbb{N}_0$ and $r>0$ ($r=1$).}
\label{tab:all_typeI_cosech_distributions}
\end{table}

\section{The Discrete Hyperbolic Cosecant Distribution -- Type I}\label{ch:I_Cosech}

We introduce the discrete type I hyperbolic cosecant distribution of the first, second, and third kind which are mainly based on the hyperbolic cosecant series identities given in \cite{ramanujan_1991} and \cite{bruckman_1977}.

We begin with the basic case, viz. the bilateral lemniscatic $\mathsf{Cosech}$ distribution. Let $w(x):=x \operatorname{cosech}\left(x \pi / 2\right)$ be the underlying continuous real weight function; we replace its removable singularity in $0$ by setting $w(0):= \lim_{x \rightarrow 0} w(x) = 2/\pi$. The discrete weight functions $w_{X_1}$, $w_{Y_1}$, and $w_{Z_1}$ are obtained by sampling the continuous weight function $w$ on the corresponding even, odd, and integer lattices, respectively. Specifically, for $\kappa\in\mathbb{Z}$, define $w_{X_1}(\kappa):=w(2\kappa)$, $w_{Y_1}(\kappa):=w(2\kappa+1)$, and $w_{Z_1}(\kappa):=w(\kappa)$ for the type I $\mathsf{Cosech}$ distributions of the first, second, and third kind, respectively.

The proofs for the normalization constants of the lemniscatic type I $\mathsf{Cosech}$ distribution are based on quasi-modularity and a duplication formula for the Eisenstein series $E_2$ for the first kind; we cite Ramanujan's formula (xiii) in Entry 15 in Chapter 17 of his second notebook for the second kind, cf. \cite[p.\,133]{ramanujan_1991}; finally, a simple parity split yields the normalization constant for the third kind. For the type I $\mathsf{Cosech}$ distribution of the first kind, the underlying weight function is $w_{X_1}(x):=w(2x)$. Its simple Fourier transform, viz.
\[
\mathcal{F}[w_{X_1}](\xi)
:=\int_{-\infty}^{\infty}2x\operatorname{cosech}(\pi x)
  e^{-2\pi i x\xi}\,dx
=\operatorname{sech}^2(\pi\xi),
\]
will be particularly useful in the context of Poisson summation. We therefore call the $\mathsf{Cosech}_{X_1}$ distribution simply the discrete hyperbolic cosecant distribution $\mathsf{Cosech}$.

\begin{propdefi}\label{prop:lemniscatic_sech}[Lemniscatic Type I $\mathsf{Cosech}$ Distributions.]
Let $\kappa \in \mathbb{Z}$. For the rows 2--4 of Table \ref{tab:all_typeI_cosech_distributions}, let $C_A$ denote the normalization constant of the listed $\mathsf{Cosech}_{A}$ distribution with weight function $w_{A}(\kappa)$. Then the associated sequence $\left(p_{A}(\kappa)\right)_{\kappa \in \mathbb{Z}} := \left(C_A \cdot w_{A}(\kappa)\right)_{\kappa \in \mathbb{Z}}$ is a probability mass function.

We name $\mathsf{Cosech}_{X_1}$, $\mathsf{Cosech}_{Y_1}$, and $\mathsf{Cosech}_{Z_1}$ the \textit{discrete lemniscatic type I hyperbolic cosecant distribution of the first, second, and third kind}, respectively. We omit adjectives when the context is clear. We call the $\mathsf{Cosech}_{X_1}$ distribution simply the \textit{discrete hyperbolic cosecant distribution} $\mathsf{Cosech}$.
\end{propdefi}

\begin{proof}
The lemniscatic normalization constants could easily be deduced from the general elliptic formulae presented in Proposition \ref{prop:r_sech}, by setting $K_r'/K_r=\sqrt{r} = 1$. But we want to examine some connections of the $\mathsf{Cosech}$ distribution to Eisenstein series and Jacobi theta functions, and further to Ramanujan's work; hence, now we give separate proofs for the lemniscatic case.\\

\textit{Ad} $C_{X_1}$: Unfortunately, the normalization constant $C_{X_1}$ is not provided by Entry 15 in Chapter 17 of Ramanujan's second notebook, cf. \cite[p.\,132--133]{ramanujan_1991}. Instead, Ramanujan states here formulae for odd moments $\ge 3$ of the $\mathsf{Cosech}$ distribution. We calculate the normalization constant $C_{X_1}$ as follows: First, we express $w_{X_1}(n)$ in terms of two Eisenstein series, viz. $E_2(i)$ and $E_2(i/2)$. Due to the quasi-modularity of $E_2$, we can directly evaluate $E_2(i)=3/\pi$ and write $E_2(i/2)$ in terms of $E_2(2i)$. Finally, we determine the value of $E_2(2i)$ by showing a duplication formula. This includes Jacobi theta functions whose values we can directly compute.

In detail: For every $n \in \mathbb{N}$, by definition we have
$$
\sinh(n \pi)=\frac{e^{\pi n}-e^{-\pi n}}{2} = \frac{1-e^{-2\pi n}}{2 e^{-\pi n}}.
$$
Using the identity
\begin{equation}\label{eq:x_1+x}
\frac{x}{1+x} = \frac{x}{1-x} - \frac{2x^2}{1-x^2},
\end{equation}
we can write
\begin{equation}\label{eq:lambert_cosech}
\sum_{n=1}^\infty 2n\operatorname{cosech}(n \pi) = 4 \sum_{n=1}^{\infty}
\frac{n e^{-\pi n}}{1-e^{-2\pi n}} = 4 \left(\sum_{n=1}^{\infty} \frac{n e^{-\pi n}}{1-e^{-\pi n}}
- \sum_{n=1}^{\infty} \frac{n e^{-2\pi n}}{1-e^{-2\pi n}} \right).
\end{equation}

Let $\tau\in\mathbb{H} =\{z\in\mathbb{C}:\Im(z)>0\}$ lie in the upper half-plane and define the nome $q_\tau = e^{i \pi \tau}$. The Eisenstein series $E_2$ has the Lambert-series expansion
\begin{equation*}\label{eisenstein_series_E2}
E_2\left(\tau\right) = 1-24\sum_{n=1}^{\infty} \frac{n q_\tau^{2n}}{1-q_\tau^{2n}}.
\end{equation*}
Hence,
\begin{equation*}
\sum_{n=1}^\infty 2n\operatorname{cosech}(n \pi)
= \frac{1}{6}
\left( E_2\left(i\right)-E_2\left(\frac{i}{2}\right) \right).
\end{equation*}

Let $G_2(\tau)=2\zeta(2)E_2(\tau)$ be a scaled version of $E_2$ by the factor $\zeta(2)=\pi^2/6$ of the Riemann zeta function. The quasi-modular transformation formula
$$
G_2\left(-\frac{1}{\tau}\right)
= \tau^2G_2(\tau)-2\pi i\tau
\quad \Leftrightarrow \quad
E_2\left(-\frac{1}{\tau}\right)
= \tau^2E_2(\tau)+\frac{6\tau}{\pi i}
$$
is given by Equation $(50)$ on p. 69 in \cite{apostol_1990}. This formula directly yields the values $E_2(i)=3/\pi$, and $E_2(i/2) = 12/\pi-4E_2(2i)$. Hence, we get
$$
\sum_{n=1}^\infty 2n \operatorname{cosech}(n \pi)
=
\frac{2}{3}E_2(2i) - \frac{3}{2\pi}.
$$

Now, we represent the Eisenstein series $E_2(2\tau)$ by a duplication formula which includes $E_2(\tau)$ and the Jacobi theta functions $\vartheta_3$ and $\vartheta_4$. Those are defined by
$$
\vartheta_3\left(q_\tau\right) = \sum_{k\in\mathbb Z}q_\tau^{k^2},
\quad
\vartheta_4\left(q_\tau\right) = \sum_{k\in\mathbb Z}(-1)^k q_\tau^{k^2},
$$
cf. \cite[p.~464]{whittaker_1920}. The duplication formula reads as follows.

\begin{lemma}[Duplication formula]
Let $\tau\in\mathbb{H} =\{z\in\mathbb{C}:\Im(z)>0\}$ and define the nome $q_\tau = e^{i \pi \tau}$. Then
\begin{equation*}
E_2(2\tau) = \frac{1}{2}E_2(\tau) + \frac{1}{4}
\left(\vartheta_3^4\left(q_\tau\right)+\vartheta_4^4\left(q_\tau\right) \right)
\end{equation*}
holds.
\end{lemma}

\begin{proof}
According to \cite[Proposition~3.2, pp. 70--71]{borwein_1987}, the fourth powers of the Jacobi theta functions admit the Lambert-series expansions
\begin{equation*}
\vartheta_3^4\left(q_\tau\right) = 1 + 8\sum_{n=1}^{\infty} \frac{2nq_\tau^{2n}}{1+q_\tau^{2n}} + 8\sum_{n=1}^{\infty} \frac{(2n-1)q_\tau^{2n-1}}{1-q_\tau^{2n-1}}
\end{equation*}
and
\begin{equation*}
\vartheta_4^4\left(q_\tau\right) = 1 + 8\sum_{n=1}^{\infty} \frac{2n q_\tau^{2n}}{1+q_\tau^{2n}} - 8\sum_{n=1}^{\infty} \frac{(2n-1)q_\tau^{2n-1}}{1+q_\tau^{2n-1}}.
\end{equation*}
Adding both identities yields
\begin{equation*}
\vartheta_3^4\left(q_\tau\right) + \vartheta_4^4\left(q_\tau\right) = 2 + 32\sum_{n=1}^{\infty} \frac{nq_\tau^{2n}}{1+q_\tau^{2n}} + 16\sum_{n=1}^{\infty} \frac{(2n-1)q_\tau^{2(2n-1)}}{1-q_\tau^{2(2n-1)}}.
\end{equation*}
Applying equation (\ref{eq:x_1+x}) to the first series and splitting the second series into odd and even indices gives
\begin{equation*}
\vartheta_3^4\left(q_\tau\right) + \vartheta_4^4\left(q_\tau\right) = 2 + 48\sum_{n=1}^{\infty} \frac{nq_\tau^{2n}}{1-q_\tau^{2n}} - 96\sum_{n=1}^{\infty} \frac{nq_\tau^{4n}}{1-q_\tau^{4n}}.
\end{equation*}
Since $E_2(\tau) = 1 - 24\sum_{n=1}^{\infty} \frac{nq_\tau^{2n}}{1-q_\tau^{2n}}$ and $E_2(2\tau) = 1 - 24\sum_{n=1}^{\infty} \frac{nq_\tau^{4n}}{1-q_\tau^{4n}}$, we conclude
$$
\vartheta_3^4\left(q_\tau\right)+\vartheta_4^4\left(q_\tau\right) = 4E_2(2\tau)-2 E_2(\tau),
$$
and the claim follows.
\end{proof}

Now, we use the well-known identities $\vartheta_3^2(q_\tau) = 2K(k)/\pi$ and $\vartheta_4^2(q_\tau) = 2 k'K(k)/\pi$, e.\,g. see Equation (5) in \cite{zucker_1984}. In the lemniscatic case, where $\lambda^*(1) = k_1=1/\sqrt{2}$ and $K(1/\sqrt{2}) = \varpi/\sqrt{2}$, this results in $\vartheta_3^4\left(q_i\right) = 2 G^2$ and
$\vartheta_4^4\left(q_i\right) = G^2$. Together with the Eisenstein value $E_2(i) = 3/\pi$ from above, we get $E_2(2i) = 3/(2\pi) + 3G^2/4$. Hence,
$$
\sum_{n=1}^{\infty} 2n \operatorname{cosech}(n \pi) = \frac{\pi G^2 - 1}{2\pi}.
$$
Since by definition, we have $w_{X_1}(0) = 2/\pi$, we finally deduce
\begin{equation}\label{eq:sum_cosech}
\sum_{\kappa \in \mathbb{Z}} 2\kappa \operatorname{cosech}(\kappa \pi) = \frac{2}{\pi} + 2 \sum_{n=1}^\infty 2n \operatorname{cosech}(n \pi) = G^2 + \frac{1}{\pi}
\quad \Leftrightarrow \quad C_{X_1} = \frac{\pi}{\pi G^2 + 1},
\end{equation}
as claimed.\\

\textit{Ad} $C_{Y_1}$: In the lemniscatic case, we translate Ramanujan's formulae in Ramanujan's Entry 15 in Chapter 17 of his second notebook, cf. \cite[p.\,132--133]{ramanujan_1991} into the common language via $x=1/2$, $y=\pi$, and $z=\sqrt{2}G$. Now, Equation (xiii) in this Entry reads as
$$
\sum_{n=0}^\infty \frac{(2n + 1)}{\sinh\left(\frac{1}{2}(2n +1)y\right)} = \frac{1}{2} z^2 \sqrt{x} = \sum_{n=0}^\infty (2n+1) \operatorname{cosech}\left((2n+1) \frac{\pi}{2} \right) = \frac{G^2}{\sqrt{2}}.
$$
Finally,
$$
\sum_{\kappa \in \mathbb{Z}} (2\kappa-1) \operatorname{cosech}\left((2\kappa-1)\frac{\pi}{2}\right) = 2\sum_{n=0}^\infty (2n+1) \operatorname{cosech}\left((2n+1) \frac{\pi}{2} \right) = \sqrt{2}G^2
$$
yields the desired normalization constant $C_{Y_1}=1/\left(\sqrt{2}G^2\right)$.\\

\textit{Ad} $C_{Z_1}$: Pairing even and odd integers yields
\begin{align*}
    \sum_{\kappa \in \mathbb{Z}} \kappa \operatorname{cosech}\left(\kappa \frac{\pi}{2} \right) &= \sum_{\kappa \in \mathbb{Z}} 2\kappa \operatorname{cosech}\left(2 \kappa \frac{\pi}{2}\right)
+
\sum_{\kappa \in \mathbb{Z}} (2\kappa - 1) \operatorname{cosech}\left((2 \kappa - 1) \frac{\pi}{2}\right)\\
&= \left(1 + \sqrt{2}\right)G^2 + \frac{1}{\pi},
\end{align*}
due to absolute convergence of each series. Hence, we get
$$C_{Z_1} = \frac{\pi}{\left(1 + \sqrt{2} \right)\pi G^2 + 1}$$
as stated.
\end{proof}

Next, we observe the discrete type I $\mathsf{Cosech}$ distribution in the general elliptic setting.

\begin{propdefi}\label{prop:r_sech}[Type I $\mathsf{Cosech}$ Distributions.]
Let $\kappa \in \mathbb{Z}$ and $r>0$. For the rows 5--7 in Table \ref{tab:all_typeI_cosech_distributions}, let $C_A$ denote the normalization constant of the listed $\mathsf{Cosech}_{A}$ distribution with weight function $w_{A}(\kappa)$. Then the associated sequence $\left(p_{A}(\kappa)\right)_{\kappa \in \mathbb{Z}} := \left(C_A \cdot w_{A}(\kappa)\right)_{\kappa \in \mathbb{Z}}$ is a probability mass function.

We name $\mathsf{Cosech}_{X_r}$, $\mathsf{Cosech}_{Y_r}$, and $\mathsf{Cosech}_{Z_r}$ the \textit{discrete type I hyperbolic cosecant distribution of the first, second, and third kind}, respectively. We omit adjectives when the context is clear.
\end{propdefi}

\begin{proof} \textit{Ad} $C_{X_r}$: For the calculation of the normalization constant $C_{X_r}$, we follow the method in \cite{bruckman_1977}. We use the Jacobi elliptic function $\operatorname{dn}$ and the Lambert-series expansion of $\operatorname{dn}^2$ for the nome $e^{-\pi K'/K}$, viz.
$$
\operatorname{dn}^2(u,k) = \frac{E}{K} + \frac{2\pi^2}{K^2} \sum_{n=1}^\infty \frac{n q^n}{1-q^{2n}} \cos\left(\frac{u}{K} \pi n\right),
$$
cf. Equation (14) in \cite[p.\,294]{bruckman_1977}. Since $\operatorname{dn}(0,k)=1$, a modification of equation (\ref{eq:lambert_cosech}) yields the identities
\begin{equation}\label{eq:n_cosech_K'_K}
4 \sum_{n=1}^\infty \frac{n q^n}{1-q^{2n}} = \sum_{n=1}^\infty 2n \operatorname{cosech}\left(\frac{K'}{K} n \pi \right) = \frac{2K(K-E)}{\pi^2},
\end{equation}
cf. Equation (25) in \cite[p.\,294]{bruckman_1977}. We choose the elliptic lambda-star function $\lambda^*(r) = k_r$ such that $K(k'_r)/K(k_r) = \sqrt{r}$. By definition of the weight function in $0$ and due to Equation \eqref{eq:E_r_alpha}, we get
\begin{align*}
\sum_{\kappa \in \mathbb{Z}} 2\kappa \operatorname{cosech}\left(\sqrt{r} \kappa \pi \right) &= \frac{2}{\sqrt{r}\pi} + 2 \sum_{n=1}^\infty 2n \operatorname{cosech}\left(\sqrt{r} n \pi \right)\\
&= \frac{2 \pi + 4 \sqrt{r} K_r (K_r - E_r)}{\sqrt{r} \pi^2} = \frac{\pi + 4 \alpha(r) K^2_r}{\sqrt{r} \pi^2}.
\end{align*}
This provides the desired normalization constant
$$
C_{X_r} = \frac{\sqrt{r} \pi^2}{\pi + 4 \alpha(r) K^2_r}.
$$

\textit{Ad} $C_{Y_r}$: Equation (xiii) of Ramanujan's Entry 15 in Chapter 17 of his second notebook, cf. \cite[p.\,133]{ramanujan_1991}, can be translated into the common language via $x=k^2_r$, $y=\pi \sqrt{r}$, and $z=2 K_r / \pi$, cf. \cite[p.\,102]{ramanujan_1991}, viz.
$$
\sum_{\kappa=0}^\infty \frac{(2\kappa + 1)}{\sinh\left(\frac{1}{2}(2\kappa +1)y\right)} = \frac{1}{2} z^2 \sqrt{x} = \sum_{n=0}^\infty (2n+1) \operatorname{cosech}\left(\sqrt{r} (2n+1) \frac{\pi}{2} \right) = \frac{2k_r K^2_r}{\pi^2}.
$$
Now,
$$
\sum_{\kappa \in \mathbb{Z}} (2\kappa+1) \operatorname{cosech}\left(\sqrt{r} (2\kappa+1)\frac{\pi}{2}\right) = \frac{4k_r K^2_r}{\pi^2}
$$
yields the normalization constant
$$
C_{Y_r}=\frac{\pi^2}{4k_rK^2_r},
$$
as claimed.\\

\textit{Ad} $C_{Z_r}$: Pairing odd and even integers yields
\begin{align*}
    \sum_{\kappa \in \mathbb{Z}} \kappa \operatorname{cosech}\left(\sqrt{r} \kappa \frac{\pi}{2} \right) &= \sum_{\kappa \in \mathbb{Z}} 2\kappa \operatorname{cosech}\left(\sqrt{r} (2 \kappa) \frac{\pi}{2}\right)
+
\sum_{\kappa \in \mathbb{Z}} (2\kappa + 1) \operatorname{cosech}\left(\sqrt{r}(2 \kappa + 1) \frac{\pi}{2}\right)\\
&= \frac{\pi + 4 \alpha(r) K^2_r + 4 \sqrt{r} k_r K^2_r}{\sqrt{r} \pi^2}.
\end{align*}
We deduce the normalization constant
$$
C_{Z_r} = \frac{\sqrt{r} \pi^2}{\pi + 4 \alpha(r) K^2_r + 4 \sqrt{r} k_r K^2_r},
$$
as stated.
\end{proof}

In \cite{ramanujan_1991} and \cite{bruckman_1977}, the hyperbolic series are formulated as a summation over $\mathbb{N}_0$. Hence, the proof of Proposition \ref{prop:r_sech} directly yields a unilateral version $\underline{\mathsf{Cosech}}$ of our type I $\mathsf{Cosech}$ distribution. For the lemniscatic case, we simply insert $r=1$, $\alpha(1)=1/2$, $k_1=1/\sqrt{2}$, and $K_1=\varpi/\sqrt{2}$ into the general elliptic formulae.

\begin{cordefi}\label{prop:r_sech_uni}[(Lemniscatic) Type I $\underline{\mathsf{Cosech}}$ Distribution.]
Let $n \in \mathbb{N}_0$ and $r>0$ $(r=1)$. For the rows 8--14 of Table \ref{tab:all_typeI_cosech_distributions}, let $\underline{C}_A$ denote the normalization constants of the listed $\underline{\mathsf{Cosech}}_{A}$ distribution with weight function $\underline{w}_{A}(n)$. Then the associated sequence $\left(p_{A}(n)\right)_{n \in \mathbb{N}_0} := \left(\underline{C}_A \cdot \underline{w}_{A}(n)\right)_{n \in \mathbb{N}_0}$ is a probability mass function.

We name $\underline{\mathsf{Cosech}}_{X_r}$, $\underline{\mathsf{Cosech}}_{Y_r}$, and $\underline{\mathsf{Cosech}}_{Z_r}$ the \textit{discrete unilateral (lemniscatic for $r=1$) type I hyperbolic cosecant distribution of the first, second, and third kind}, respectively. We omit adjectives when the context is clear.
\end{cordefi}

\section{The Discrete Hyperbolic Cosecant Distribution -- Type II}\label{ch:II_Cosech}

Now, we introduce the discrete type II hyperbolic cosecant distribution, building on the foundational work of \cite{ramanujan_1991}, \cite{bruckman_1977}, and \cite{zucker_1984}.

For the type I $\mathsf{Cosech}$ distribution, the underlying continuous weight function $w(x) = x \operatorname{cosech}(x\pi/2)$ has a removable singularity at $0$. The $x$ factor is necessary to make $w$ positive and even on its domain $\mathbb{Z}$. In contrast, for the type II $\mathsf{Cosech}$ distribution, we use the underlying continuous real weight function $w_{\mathsf{II}}(x) := \operatorname{cosech}(x \pi /2) / x$ which has a non-removable double pole at $0$. Hence, we first consider the unilateral type II $\underline{\mathsf{Cosech}}$ distributions on $\mathbb{N}_0$ after index shift. We prove the general elliptic formulae and for the lemniscatic case where $r=1$, we insert the elliptic singular modulus $k_1=1/\sqrt{2}$ into the formulae.

\renewcommand{\arraystretch}{2}
\begin{table}[!tbh]
    \centering
    \begin{tabular}{|c|c|c|c|}
\hline
Type & Distribution & Normalizing Constant & Weight Function\\
\hline
II & $\mathsf{\underline{Cosech}}_{U_1}$ & $\frac{24}{\pi - \ln(8)}$ & $\operatorname{cosech}\left((n+1) \pi \right) / (2n+2)$\\
II & $\mathsf{\underline{Cosech}}_{V_1}$ & $\frac{2}{\ln(1+\sqrt{2})} = \frac{2}{\operatorname{arsinh}(1)}$ & $\operatorname{cosech} \left((2n + 1) \frac{\pi}{2}\right) / (2n + 1)$\\
II & $\mathsf{\underline{Cosech}}_{W_1}$ & $\frac{24}{\pi - \ln(8) + 12\ln(1+\sqrt{2})}$ & $ \operatorname{cosech} \left((n+1) \frac{\pi}{2} \right) / (n+1)$\\
\hline
II & $\mathsf{\underline{Cosech}}_{U_r}$ & $\frac{24}{\sqrt{r}\pi + 2 \ln(k_r) - 2\ln(2 + 2k_r) - 2\ln(2 - 2k_r)}$ & $\operatorname{cosech} \left(\sqrt{r} (n+1) \pi \right) / (2n+2)$\\
II & $\mathsf{\underline{Cosech}}_{V_r}$ & $\frac{2}{\operatorname{artanh}\left(k_r\right)}$ & $\operatorname{cosech} \left(\sqrt{r} (2n + 1) \frac{\pi}{2}\right) / (2n + 1)$\\
II & $\mathsf{\underline{Cosech}}_{W_r}$ & $\frac{24}{\sqrt{r}\pi + 2\ln(4k_r) + 4\ln(1+k_r) - 8\ln (2-2k_r)}$ & $ \operatorname{cosech} \left(\sqrt{r} (n+1) \frac{\pi}{2} \right) / (n+1)$\\
\hline
II & $\mathsf{Cosech}_{V_1}$ & $\frac{1}{\operatorname{arsinh}(1)}$ & $\operatorname{cosech} \left((2\kappa + 1) \frac{\pi}{2}\right) / (2\kappa + 1)$\\
II & $\mathsf{Cosech}_{V_r}$ & $\frac{2}{\ln \left(1+k_r\right) - \ln(1-k_r)}$ & $\operatorname{cosech} \left(\sqrt{r} (2\kappa + 1) \frac{\pi}{2}\right) / (2\kappa + 1)$\\
\hline
\end{tabular}
\caption{Discrete unilateral (lemniscatic) type II $\underline{\mathsf{Cosech}}$ distributions of the first ($U$), second ($V$), and third ($W$) kind for $n \in \mathbb{N}_0$ and $r>0$ ($r=1$); further the discrete bilateral (lemniscatic) type II $\mathsf{Cosech}$ distributions of the second ($V$) kind for $\kappa \in \mathbb{Z}$ and $r>0$ ($r=1$).}
\label{tab:all_typeII_cosech_distributions}
\end{table}

\begin{propdefi}\label{prop:r_II_cosech}[(Lemniscatic) Type II $\underline{\mathsf{Cosech}}$ Distributions.]
Let $n \in \mathbb{N}_0$ and $r>0$ ($r=1$). For the rows 2--7 of Table \ref{tab:all_typeII_cosech_distributions}, let $\underline{C}_A$ denote the normalization constant of the listed $\mathsf{\underline{Cosech}}_{A}$ distribution with weight function $\underline{w}_{A}(n)$. Then the associated sequence $\left(p_{A}(n)\right)_{n \in \mathbb{N}_0} := \left(\underline{C}_A \cdot \underline{w}_{A}(n)\right)_{n \in \mathbb{N}_0}$ is a probability mass function.

We name $\mathsf{\underline{Cosech}}_{U_r}$, $\mathsf{\underline{Cosech}}_{V_r}$, and $\mathsf{\underline{Cosech}}_{W_r}$ the \textit{discrete (lemniscatic for $r=1$) unilateral type II hyperbolic cosecant distribution of the first, second, and third kind}, respectively. We omit adjectives when the context is clear.
\end{propdefi}

\begin{proof}
    \textit{Ad} $\underline{C}_{U_r}$: We follow the technique in \cite{zucker_1984}. Define the nome $q=e^{- \pi K'/K}$ with the complete elliptic integral of the first kind $K$ and further the infinite products
$$
Q_0 := \prod_{n=1}^\infty \left(1-q^{2n}\right), \quad Q_3 := \prod_{n=1}^\infty \left(1-q^{2n-1}\right).
$$
Equation (5) in \cite[p.\,406]{zucker_1984} yields the equations
$$
Q_0^2 Q_3^4= \frac{2 k' K}{\pi}, \quad 2 q^{1/4} Q_0^3 = \sqrt{\frac{8K^3 k k'}{\pi^3}}.
$$
On the one hand, from these identities we deduce
$$
Q_3 = q^{1/24} \frac{(2 k')^{1/6}}{k^{1/12}} \quad \Leftrightarrow \quad -\ln\left(Q_3\right) = \frac{\frac{K'}{K} \pi + 2\ln(k) - 4 \ln(2 k')}{24}.
$$
On the other hand, the logarithmic series $\ln(1-x) = - \sum_{m=1}^\infty x^m/m$ and the geometric series results in
\begin{align*}
-\ln\left(Q_3\right) &= -\sum_{m=1}^\infty \ln \left(1 - q^{2m-1}\right) = \sum_{m=1}^\infty \sum_{n=1}^\infty \frac{q^{(2m-1)n}}{n} = \sum_{n=1}^\infty \frac{q^n}{n\left(1-q^{2n}\right)}\\
&= \frac{1}{2} \sum_{n=1}^\infty \frac{2}{n\left(e^{\pi n \frac{K'}{K}} - e^{-\pi n \frac{K'}{K}}\right)} = \sum_{n=1}^\infty \frac{\operatorname{cosech}\left(\frac{K'}{K}n\pi\right)}{2n}.
\end{align*}
Hence,
\begin{equation}\label{eq:cosech_n_pi}
\sum_{n=0}^\infty \frac{\operatorname{cosech}\left(\frac{K'}{K}(2n+2)\frac{\pi}{2} \right)}{2n+2} = \frac{\frac{K'}{K}\pi + 2 \ln(k)- 4 \ln(2k')}{24}.
\end{equation}

We choose the elliptic lambda-star function $\lambda^*(r) = k_r$ such that $K'_r/K_r = \sqrt{r}$. This yields
$$
\sum_{n=0}^\infty \frac{\operatorname{cosech}(\sqrt{r}(n+1)\pi)}{2n+2} = \frac{\sqrt{r}\pi + 2 \ln(k_r) - 4 \ln(2k'_r)}{24}.
$$
This results in the normalization constant
$$
\underline{C}_{U_r} = \frac{24}{\sqrt{r}\pi + 2 \ln(k_r) - 2\ln(2 + 2k_r) - 2\ln(2 - 2k_r)},
$$
as claimed.\\

\textit{Ad} $\underline{C}_{V_r}$: The suitable formula is given by the Example after Entry 7 of Chapter 18 in Ramanujan's second notebook, cf. \cite[p.\,156]{ramanujan_1991}. With $x=k^2$ and $y=\pi K'/K$, Ramanujan states
\begin{align*}
& \sum_{n=0}^\infty \frac{1}{(2n+1) \sinh \left(\frac{1}{2}(2n+1) y\right)} = \frac{1}{4} \log \left(\frac{1 + \sqrt{x}}{1 - \sqrt{x}}\right)\\
=& \sum_{n=0}^\infty \frac{\operatorname{cosech}\left(\frac{K'}{K}(2n+1)\frac{\pi}{2}\right)}{2n+1} = \frac{1}{4} \ln \left(\frac{1 + k}{1 - k}\right).
\end{align*}
This is e.\,g. proven by Berndt in \cite[p.\,156--157]{ramanujan_1991}. We use Equation (17) in \cite[p.\,294]{bruckman_1977} with $k \in (0,1)$ and $q = e^{-\pi K' / K}$, viz.
$$
-\frac{1}{16}\ln(1-k^2)
=
\sum_{n=0}^{\infty}
\frac{q^{2n+1}}{(2n+1)(1-q^{4n+2})}.
$$
For each $n\in\mathbb{N}_0$,
$$
\operatorname{cosech}\left(\frac{K'}{K}(2n+1)\pi\right)
=
\frac{2q^{2n+1}}{1-q^{4n+2}}.
$$
Thus
$$
\begin{aligned}
\sum_{n=0}^{\infty}
\frac{\operatorname{cosech}\left(\frac{K'}{K}(2n+1)\pi\right)}{2n+1} = 2\sum_{n=0}^{\infty}
\frac{q^{2n+1}}{(2n+1)(1-q^{4n+2})} = -\frac{1}{4}\ln(k').
\end{aligned}
$$
The series converges because \(0<q<1\). In \cite[p.\,409]{zucker_1984}, the following transformations for changing the nome from $q$ to $q^{1/2}$ are stated:
\begin{equation}\label{eq:zucker_trafo}
k \rightarrow \frac{2\sqrt{k}}{1+k}, \quad k' \rightarrow \frac{1-k}{1+k}, \quad \frac{K'}{K} \rightarrow \frac{K'}{2K}.
\end{equation}
We again choose the elliptic lambda-star function $\lambda^*(r) = k_r$ such that $K'_r/K_r = \sqrt{r}$. This yields
$$
\sum_{n=0}^\infty \frac{\operatorname{cosech}\left(\sqrt{r}(2n+1)\frac{\pi}{2}\right)}{2n+1} = \frac{1}{4} \ln \left(\frac{1+k_r}{1-k_r}\right) = \frac{\operatorname{artanh}(k_r)}{2}.
$$
Therefore, we get the normalization constant
$$
\underline{C}_{V_r} = \frac{4}{\ln\left(1+k_r\right) - \ln\left(1-k_r\right)} = \frac{2}{\operatorname{artanh}(k_r)},
$$
as claimed.\\

\textit{Ad} $\underline{C}_{W_r}$: We use the transformations stated in Equation \eqref{eq:zucker_trafo} from \cite[p.\,409]{zucker_1984}. We choose the elliptic lambda-star function $\lambda^*(r) = k_r$ such that $K'_r/K_r = \sqrt{r}$. Hence, Equation \eqref{eq:cosech_n_pi} yields
$$
\sum_{n=0}^\infty \frac{\operatorname{cosech}\left(\sqrt{r}(n + 1)\frac{\pi}{2}\right)}{n+1} = \frac{\sqrt{r}\pi + 2\ln(4k_r) + 4\ln(1+k_r) - 8\ln(2-2k_r)}{24}.
$$
Thus, we finally get the normalization constant
$$
\underline{C}_{W_r} = \frac{24}{\sqrt{r}\pi + 2\ln(4k_r) + 4\ln(1+k_r) - 8\ln(2-2k_r)}
$$
as stated.
\end{proof}

Since it leaves the $0$ out, we can also define a bilateral type II $\mathsf{Cosech}$ of the second kind. By symmetry, we just have to take the normalization constant $C_{V_r} = \underline{C}_{V_r}/2$.

\begin{cordefi}\label{cor:bilateral_II_cosech}[Type II $\mathsf{Cosech}$ Distributions.]
Let $\kappa \in \mathbb{Z}$. For the rows 8--9 of Table \ref{tab:all_typeII_cosech_distributions}, let $C_A$ denote the normalization constant of the listed $\mathsf{Cosech}_{A}$ distribution with weight function $w_{A}(\kappa)$. Then the associated sequence $\left(p_A(\kappa)\right)_{\kappa \in \mathbb{Z}} := \left(C_A \cdot w_{A}(\kappa)\right)_{\kappa \in \mathbb{Z}}$ is a probability mass function.

We name $\mathsf{Cosech}_{V_r}$ and $\mathsf{Cosech}_{V_1}$ the \textit{discrete bilateral (lemniscatic for $r=1$) type II hyperbolic cosecant distribution of the second kind}, respectively. We omit adjectives when the context is clear.
\end{cordefi}

\section{The Discrete Hyperbolic Cosecant Distribution -- Type III}

Next, we introduce a further version of the $\mathsf{Cosech}$ distribution: The type III hyperbolic cosecant distribution. Like for the type II $\mathsf{Cosech}$ distribution, the underlying continuous weight function $w_\mathsf{III}(x) := \operatorname{cosech}^2(x \pi)$ has a non-removable double pole at $0$. Therefore, we again first concentrate on the unilateral type III $\underline{\mathsf{Cosech}}$ distributions on $\mathbb{N}_0$ after index shift. We use techniques from \cite{ramanujan_1991}, \cite{borwein_1987}, and \cite{zucker_1979} to prove the general elliptic case for $r>0$. For the lemniscatic case, we insert $r=1$, $k_1=1/\sqrt{2}$, $K_1 = \varpi/\sqrt{2}$, and $\alpha(1)=1/2$ into the general formulae.

\renewcommand{\arraystretch}{2}
\begin{table}[!tbh]
    \centering
    \begin{tabular}{|c|c|c|c|}
\hline
Type & Distribution & Normalizing Constant & Weight Function\\
\hline
III & $\mathsf{\underline{Cosech}}_{R_1}$ & $\frac{6 \pi}{\pi -3}$ & $\operatorname{cosech}^2\left((2n+2) \frac{\pi}{2}\right)$\\
III & $\mathsf{\underline{Cosech}}_{S_1}$ & $\frac{2 \pi}{\pi G^2 - 1}$ & $\operatorname{cosech}^2 \left((2n + 1) \frac{\pi}{2}\right)$\\
III & $\mathsf{\underline{Cosech}}_{T_1}$ & $\frac{6 \pi}{\pi - 6 + 3 \pi G^2}$ & $\operatorname{cosech}^2\left((n+1) \frac{\pi}{2}\right)$\\
\hline
III & $\mathsf{\underline{Cosech}}_{R_r}$ & $\frac{6 \sqrt{r} \pi^2}{\sqrt{r}\pi^2 -3 \pi + 4K_r^2\left(3 \alpha(r) - (1 + k_r^2)\sqrt{r}\right)}$ & $\operatorname{cosech}^2\left(\sqrt{r} (2n+2) \frac{\pi}{2}\right)$\\
III & $\mathsf{\underline{Cosech}}_{S_r}$ & $\frac{2\sqrt{r}\pi^2}{4 K_r^2 \alpha(r) - \pi}$ & $\operatorname{cosech}^2 \left(\sqrt{r} (2n + 1) \frac{\pi}{2}\right)$\\
III & $\mathsf{\underline{Cosech}}_{T_r}$ & $\frac{6\sqrt{r}\pi^2}{\sqrt{r}\pi^2- 6\pi + 4 K_r^2 \left(6 \alpha(r) - (1 + k_r^2) \sqrt{r}\right)}$ & $\operatorname{cosech}^2\left(\sqrt{r} (n+1) \frac{\pi}{2}\right)$\\
\hline
III & $\mathsf{Cosech}_{S_1}$ & $\frac{\pi}{\pi G^2 - 1}$ & $\operatorname{cosech}^2 \left((2\kappa + 1) \frac{\pi}{2}\right)$\\
III & $\mathsf{Cosech}_{S_r}$ & $\frac{\sqrt{r}\pi^2}{4 K_r^2 \alpha(r) - \pi}$ & $\operatorname{cosech}^2 \left(\sqrt{r} (2\kappa + 1) \frac{\pi}{2}\right)$\\
\hline
\end{tabular}
\caption{Discrete unilateral (lemniscatic) type III $\underline{\mathsf{Cosech}}$ distributions of the first ($R$), second ($S$), and third ($T$) kind for $n \in \mathbb{N}$ and $r>0$ ($r=1$); further the discrete bilateral (lemniscatic) type III $\mathsf{Cosech}$ distributions of the second ($S$) kind for $\kappa \in \mathbb{Z}$ and $r>0$ ($r=1$).}
\label{tab:all_typeIII_cosech_distributions}
\end{table}

\begin{propdefi}\label{prop:r_III_uni_cosech}[(Lemniscatic) Type III $\underline{\mathsf{Cosech}}$ Distributions.]
Let $n \in \mathbb{N}_0$ and $r>0$ ($r=1$). For the rows 2--6 of Table \ref{tab:all_typeIII_cosech_distributions}, let $\underline{C}_A$ denote the normalization constant of the listed $\mathsf{\underline{Cosech}}_{A}$ distribution with weight function $\underline{w}_{A}(n)$. Then the associated sequence $\left(p_{A}(n)\right)_{n \in \mathbb{N}_0} := \left(\underline{C}_A \cdot \underline{w}_{A}(n)\right)_{n \in \mathbb{N}_0}$ is a probability mass function.

We name $\mathsf{\underline{Cosech}}_{R_r}$, $\mathsf{\underline{Cosech}}_{S_r}$, and $\mathsf{\underline{Cosech}}_{T_r}$ the \textit{discrete (lemniscatic for $r=1$) unilateral type III hyperbolic cosecant distribution of the first, second, and third kind}, respectively. We omit adjectives when the context is clear.
\end{propdefi}

\begin{proof}
\textit{Ad} $\underline{C}_{R_r}$: For the calculation of the normalizing constant $\underline{C}_{R_r}$ for the $\mathsf{\underline{Cosech}}_{R_r}$ distribution we can use the technique in \cite[p.\,197]{zucker_1979} or complete Exercise $7$ of chapter 3.7 in \cite[p.\,96]{borwein_1987}. We do the latter. Let $c:=K'_r/K_r$ with nome $q=e^{-\pi c}$, where $k=k(c)$ is determined implicitly. Further, define
$$
Q_0:=\prod_{n=1}^{\infty}(1-q^{2n}).
$$
Equation $(3.2.28i)$ in \cite[p.\,74]{borwein_1987} claims
$$
F(c)
:= -\ln\left(Q_0\right)
=\frac12\sum_{n=1}^{\infty}\frac{\coth(c n \pi)-1}{n}
=-\frac{c \pi}{12} -\frac{1}{6}\ln\left(\frac{2kk'K^3}{\pi^3}\right).
$$
For \(c>0\), the series and its termwise derivative converge uniformly on every compact subinterval of \((0,\infty)\). Hence, the derivative of the series yields
$$
F'(c)
=-\frac{\pi}{2}\sum_{n=1}^{\infty}\operatorname{cosech}^2(c n \pi).
$$
Differentiating the right-hand side with respect to $c$ results in
$$
F'(c)
=-\frac{\pi}{12}
-\frac{1}{6}\frac{\mathrm{d}}{\mathrm{d}c}\ln(2kk'K^3).
$$
Equating the two expressions for \(F'(c)\), we obtain
\begin{equation*}
\sum_{n=1}^{\infty}\operatorname{cosech}^2(c n \pi)
=
\frac16+\frac{1}{3\pi}\frac{\mathrm{d}}{\mathrm{d}c}\ln(2kk'K^3).
\end{equation*}
We now compute the derivative on the right-hand side. Since \(K\) is a function of \(k\), and \(k\) is a function of \(c\), the expression \(f(K,k)\) depends on \(c\) through \(k(c)\). Thus, by the chain rule,
$$
\frac{\mathrm{d}}{\mathrm{d}c}f(K(k(c)),k(c))
=
\left(
\frac{\partial f}{\partial K}\frac{\mathrm{d}K}{\mathrm{d}k}
+
\frac{\partial f}{\partial k}
\right)\frac{\mathrm{d}k}{\mathrm{d}c}.
$$
The expression in parentheses is the total derivative with respect to $k$. Hence,
$$
\frac{\mathrm{d}}{\mathrm{d}c}\ln(2kk'K^3)
=
\frac{\mathrm{d}}{\mathrm{d}k}\ln(2kk'K^3)\frac{\mathrm{d}k}{\mathrm{d}c}.
$$
The standard derivative formula, given in equation 19.4.1 of \cite{nist_2026}, is
$$
\frac{\mathrm{d}K}{\mathrm{d}k} = \frac{E-k'^2K}{kk'^2}.
$$
This implies
\begin{equation*}
\frac{\mathrm{d}}{\mathrm{d}k}\ln(2kk'K^3)
=3\frac{1}{K}\frac{\mathrm{d}K}{\mathrm{d}k}
  +\frac{1}{k}+\frac{1}{k'}\frac{\mathrm{d}k'}{\mathrm{d}k}
=\frac{3E}{kk'^2K}-\frac{2}{k}-\frac{k}{k'^2}.
\end{equation*}

Since $c=K(k')/K(k)$, differentiating this quotient with respect to $k$ and using the derivative formula for $K$ gives
$$
\frac{\mathrm{d}c}{\mathrm{d}k}
= -\frac{K E'+K'E-KK'}{kk'^2K^2}
= -\frac{\pi}{2kk'^2K^2},
$$
where the last equality follows from Legendre's relation
$$
EK'+E'K-KK'=\frac{\pi}{2}.
$$
Therefore, inversion yields
$$
\frac{\mathrm{d}k}{\mathrm{d}c}
=-\frac{2kk'^2K^2}{\pi}.
$$
Altogether,
$$
\frac{\mathrm{d}}{\mathrm{d}c}\ln(2kk'K^3)
=\left(\frac{3E}{kk'^2K}-\frac2k-\frac{k}{k'^2}\right)
  \left(-\frac{2kk'^2K^2}{\pi}\right)
= \frac{(4-2k^2)K^2 - 6KE}{\pi}.
$$
We choose again the elliptic lambda-star function $\lambda^*(r) = k_r$ such that $c= K'_r/K_r = \sqrt{r}$. Next, we use the formula
$$
E_r = \frac{\pi}{4\sqrt{r}K_r} + \left(1-\frac{\alpha(r)}{\sqrt{r}}\right)K_r
$$
with the elliptic alpha function $\alpha(r)$. This results in
\begin{align*}
\sum_{n=0}^{\infty}\operatorname{cosech}^2\left(\sqrt{r} (2n+2) \frac{\pi}{2}\right)
&=
\frac{\pi^2-12K_rE_r+(8-4k_r^2)K_r^2}{6\pi^2}\\
&= \frac{\sqrt{r}\pi^2 -3 \pi + 4K_r^2\left(3 \alpha(r) - (1 + k_r^2)\sqrt{r}\right)}{6 \sqrt{r} \pi^2}.
\end{align*}

Therefore,
$$
\underline{C}_{R_r} = \frac{6 \sqrt{r} \pi^2}{\sqrt{r}\pi^2 -3 \pi + 4K_r^2\left(3 \alpha(r) - (1 + k_r^2)\sqrt{r}\right)},
$$
as claimed.\\

\textit{Ad} $\underline{C}_{S_r}$: Berndt provides the Equation (14.15) in the proof of Entry 14 in Chapter 18 of Ramanujan's second notebook, cf. \cite[p.\,172]{ramanujan_1991}, viz.
\begin{align*}
& \sum_{n=-\infty}^\infty \operatorname{cosech} \left\{\frac{1}{2}(2n+1)y \right\}  = 4 \sum_{n=1}^\infty \frac{n}{\sinh(ny)} \\
=&\ \sum_{n\in\mathbb{Z}} \operatorname{cosech}^2\left(\frac{K'}{K}(2n+1)\frac{\pi}{2}\right) = 4 \sum_{n=1}^\infty n\operatorname{cosech}\left(\frac{K'}{K} n \pi\right).
\end{align*}
By the symmetry $n \mapsto -n-1$, the left-hand side equals
$$
2\sum_{n=0}^\infty \operatorname{cosech}^2\left(\frac{K'}{K}(2n+1)\frac{\pi}{2}\right)
$$
By Equation \eqref{eq:n_cosech_K'_K}, the right-hand side equals $4K(K-E)/\pi^2$. We choose the elliptic lambda-star function $\lambda^*(r) = k_r$ such that $K'_r/K_r = \sqrt{r}$. Now, Equation \eqref{eq:E_r_alpha} yields
$$
\sum_{n=0}^\infty \operatorname{cosech}^2\left(\sqrt{r}(2n+1)\frac{\pi}{2}\right) = \frac{2K_r(K_r-E_r)}{\pi^2} = \frac{4\alpha(r)K_r^2 - \pi}{2\sqrt{r}\pi^2}
$$
with the elliptic alpha function $\alpha(r)$. Therefore, we get the normalization constant
$$
\underline{C}_{S_r} = \frac{2\sqrt{r}\pi^2}{4\alpha(r)K_r^2 - \pi},
$$
as claimed.\\

\textit{Ad} $\underline{C}_{T_r}$: We split the sum over $n \geq 1$ into even and odd indices. Then
$$
\sum_{n=0}^\infty \operatorname{cosech}^2\left(\sqrt{r} (n+1) \frac{\pi}{2}\right)
= \sum_{n=0}^\infty \operatorname{cosech}^2\left(\sqrt{r} (2n+2) \frac{\pi}{2}\right)
+ \sum_{n=0}^\infty \operatorname{cosech}^2\left(\sqrt{r}(2n+1)\frac{\pi}{2}\right).
$$
By the proofs for $\underline{C}_{R_r}$ and $\underline{C}_{S_r}$, this yields
\begin{align*}
\sum_{n=0}^\infty \operatorname{cosech}^2\left(\sqrt{r} (n+1) \frac{\pi}{2}\right) &= \frac{\sqrt{r}\pi^2 - 3\pi + 4K_r^2\left(3\alpha(r) - (1+k_r^2)\sqrt{r}\right)}{6\sqrt{r}\pi^2}
+ \frac{4\alpha(r)K_r^2 - \pi}{2\sqrt{r}\pi^2}\\
&= \frac{\sqrt{r}\pi^2 - 6\pi + 4K_r^2\left(6\alpha(r) - (1+k_r^2)\sqrt{r}\right)}{6\sqrt{r}\pi^2}.
\end{align*}
Therefore,
$$
\underline{C}_{T_r} = \frac{6\sqrt{r}\pi^2}{\sqrt{r}\pi^2 - 6\pi + 4K_r^2\left(6\alpha(r) - (1+k_r^2)\sqrt{r}\right)},
$$
as stated.
\end{proof}

Since its support excludes $0$, we can also define the bilateral type III $\mathsf{Cosech}$ of the second kind. By symmetry, we just have to take the normalization constant $C_{S_r} = \underline{C}_{S_r}/2$.

\begin{cordefi}\label{cor:bilateral_lemniscatic_III_cosech}[Type III $\mathsf{Cosech}$ Distributions.]
Let $\kappa \in \mathbb{Z}$. For the rows 8--9 of Table \ref{tab:all_typeIII_cosech_distributions}, let $C_A$ denote the normalization constant of the listed $\mathsf{Cosech}_{A}$ distribution with weight function $w_{A}(\kappa)$. Then the associated sequence $\left(p_A(\kappa)\right)_{\kappa \in \mathbb{Z}} := \left(C_A \cdot w_{A}(\kappa)\right)_{\kappa \in \mathbb{Z}}$ is a probability mass function.

We name $\mathsf{Cosech}_{S_r}$ and $\mathsf{Cosech}_{S_1}$ the \textit{discrete bilateral (lemniscatic for $r=1$) type III hyperbolic cosecant distribution of the second kind}, respectively. We omit adjectives when the context is clear.
\end{cordefi}

\section{Moments of the Discrete Hyperbolic Cosecant Distribution}

We want to determine some moments of the discrete hyperbolic cosecant distribution $\mathsf{Cosech}$. According to Table \ref{tab:all_typeI_cosech_distributions}, its probability mass function is given by
$$
\left(p(\kappa)\right)_{\kappa \in \mathbb{Z}} = \frac{\pi}{\pi G^2 + 1} 2\kappa \operatorname{cosech}(\kappa \pi).
$$
Ramanujan's series identities directly yield the even moments up to the eighth order.

\begin{prop}\label{prop:Var_cosech}
Let $X\sim\mathsf{Cosech}$ be a discrete random variable which is distributed according to the hyperbolic cosecant distribution. The odd moments vanish and we have
\begin{align*}
(i) \quad \mathbb{E}\left[X^2 \right] =&\ \operatorname{Var}(X) = \frac{\pi G^4}{\pi G^2 + 1} \approx 0.4781246,\\
(ii) \quad \mathbb{E}\left[X^4 \right] =&\ \frac{3\pi G^6}{\pi G^2 + 1}\approx 0.9991876,\\
(iii) \quad \mathbb{E}\left[X^6 \right] =&\ \frac{18 \pi G^8}{\pi G^2 + 1}\approx 4.1762161,\\
(iv) \quad \mathbb{E}\left[X^8 \right] =&\ \frac{189\pi G^{10}}{\pi G^2 + 1} \approx 30.5461835,
\end{align*}
with the Gau\ss ian constant $G = \varpi/\pi$.
\end{prop}

\begin{proof}
Let $C = \pi/(\pi G^2 + 1)$. We get $\mathbb{E}\left[X^j\right] = 0$ for odd $j \in \mathbb{N}$ due to symmetry. The other formulae directly follow by Ramanujan's Entry 15 (i)-(iv) in Chapter 17 of his second notebook, cf. \cite[pp.\,132-133]{ramanujan_1991}: For $j\in \{2,4,6,8\}$ we write
$$
\mathbb{E}\left[X^{j}\right] = C \sum_{\kappa\in\mathbb{Z}} \kappa^{j} \cdot 2\kappa \operatorname{cosech}(\kappa \pi) = \frac{4\pi}{\pi G^2 + 1} \sum_{n=1}^\infty \frac{n^{j+1}}{\sinh(\pi n)}.
$$
In the lemniscatic case, we can translate Ramanujan's notation via $x=1/2$, $y=\pi$, and $z=\sqrt{2}G$. We get
$$
\sum_{n=1}^\infty \frac{n^{3}}{\sinh(y n)} = \frac{1}{8} z^4 x
\ \Rightarrow \
\mathbb{E}\left[X^2\right] = \frac{\pi G^4}{\pi G^2 + 1},
$$
$$
\sum_{n=1}^\infty \frac{n^{5}}{\sinh(y n)} = \frac{1}{8} z^6 x(1+x)
\ \Rightarrow \
\mathbb{E}\left[X^4\right] = \frac{3\pi G^6}{\pi G^2 + 1},
$$
$$
\sum_{n=1}^\infty \frac{n^{7}}{\sinh(y n)} = \frac{1}{8} z^8 x\left(1+\frac{13}{2}x+x^2\right)
\ \Rightarrow \
\mathbb{E}\left[X^6\right] = \frac{18 \pi G^8}{\pi G^2 + 1},
$$
$$
\sum_{n=1}^\infty \frac{n^{9}}{\sinh(y n)} = \frac{1}{8} z^{10} x(1+x)\left(1 + 29x + x^2 \right)
\ \Rightarrow \
\mathbb{E}\left[X^8\right] = \frac{189\pi G^{10}}{\pi G^2 + 1},
$$
as claimed.
\end{proof}

We want to evaluate some exact expressions for series of even powers of the hyperbolic secant via the Poisson summation. For that purpose we first prove a formula for their Fourier transform.

\begin{lemma} \label{le:fourier_sech}
Let $m\in\mathbb{N}$ be a positive integer and let $\xi \in \mathbb{R}$. Then the Fourier transformation formula
$$
\mathcal{F}\left[x\mapsto \operatorname{sech}^{2m}(x \pi)\right](\xi) = \int_{-\infty}^\infty \operatorname{sech}^{2m}(x\pi) e^{-2 \pi i x \xi}\,\mathrm{d}x = \frac{2^{2m-1}}{\pi (2m-1)!} \left|\Gamma(m + i\xi)\right|^2
$$
holds for every $m$. Further, the identity
$$
\left|\Gamma(m+i\xi)\right|^2 = \xi \pi \operatorname{cosech}(\xi \pi) \prod_{j=1}^{m-1}(j^2 + \xi^2)
$$
holds.
\end{lemma}

\begin{proof}
We use the hyperbolic substitution $u = u(x) := (1 + \tanh(x\pi))/2$, and hence $2u-1 = \tanh(x\pi)$, which implies $x=\ln(u/(1-u))/(2\pi)$. Now, $4 u(1-u) = 1-\tanh^2(x \pi) = \operatorname{sech}^2(x\pi)$ leads to
$$
\operatorname{sech}^{2m}(x\pi) = 4^m u^m (1-u)^m.
$$
For the boundaries, we get $\lim_{x\rightarrow -\infty}u(x) = 0$ and $\lim_{x\rightarrow \infty}u(x) = 1$. Hence, by continuity, $u\in(0,1)$ and thus, $u/(1-u) \in (0,\infty)$. Therefore, we can use the unique, real natural logarithm to get
$$
e^{-2 \pi i x \xi} = \exp\left(-i \xi \ln\left(\frac{u}{1-u}\right) \right) = u^{-i\xi} (1-u)^{i\xi}.
$$
Further, $u'(x) = \pi \operatorname{sech}^2(x\pi)/2$ results in
$$
\mathrm{d}x = \frac{2\mathrm{d}u}{\pi \operatorname{sech}^2(x\pi)} = \frac{\mathrm{d}u}{2\pi u(1-u)}.
$$
We use Euler's beta function $B: \{(a,b) \in \mathbb{C}^2: \Re(a), \Re(b)>0\} \rightarrow \mathbb{C}$ with its key property
$$
B(a,b) = \int_0^1 t^{a-1} (1-t)^{b-1} \,\mathrm{d}t = \frac{\Gamma(a)\cdot\Gamma(b)}{\Gamma(a+b)}.
$$
Altogether, we get
\begin{align*}
\int_{-\infty}^\infty \operatorname{sech}^{2m}(x\pi) e^{-2 \pi i x \xi}\,\mathrm{d}x &= \int_0^1 \underbrace{2^{2m}u^m(1-u)^m}_{\operatorname{sech}^{2m}(x\pi)}\cdot \underbrace{u^{-i\xi}(1-u)^{i\xi}}_{e^{-2 \pi i x \xi}} \frac{\mathrm{d}u}{2\pi u(1-u)}\\
&= \frac{2^{2m-1}}{\pi} \int_0^1 u^{m -i\xi -1} (1-u)^{m+i\xi-1}\, \mathrm{d}u\\
&= \frac{2^{2m-1}}{\pi} B(m -i\xi, m +i\xi)\\
&= \frac{2^{2m-1}}{\pi} \frac{\overline{\Gamma(m+i\xi)}\cdot \Gamma(m+i\xi)}{\Gamma(2m)}\\
&= \frac{2^{2m-1}}{\pi (2m-1)!} \left|\Gamma(m + i\xi)\right|^2,
\end{align*}
as claimed.\\

The second formula follows from induction:\\

\textbf{m=1}: Euler's reflection formula yields
$$
\Gamma(i\xi)\cdot\Gamma(1-i\xi) = \frac{\pi}{\sin(i \xi \pi)} = \frac{\pi}{i \sinh(\xi \pi)}.
$$
Since $\Gamma(1-i\xi) = -i\xi\Gamma(-i\xi) = - i \xi \overline{\Gamma(i\xi)}$, we get $\left| \Gamma(i\xi) \right|^2 = \pi \operatorname{cosech}(\xi \pi)/\xi$. This results in
$$
\left|\Gamma(1+i\xi)\right|^2 = \xi^2|\Gamma(i\xi)|^2 = \xi \pi \operatorname{cosech}(\xi \pi),
$$
as claimed by defining the empty product as $1$. The case $\xi=0$ follows from continuity.\\

\textbf{m} $\rightarrow$ \textbf{m+1}: Let
$$
\left|\Gamma(m+i\xi)\right|^2 = \xi \pi \operatorname{cosech}(\xi \pi) \prod_{j=1}^{m-1}(j^2 + \xi^2)
$$
hold for $m$. Then we get
\begin{align*}
\left|\Gamma(m+1+i\xi) \right|^2 &= \left(m^2 + \xi^2\right) \xi \pi \operatorname{cosech}(\xi\pi) \prod_{j=1}^{m-1} \left(j^2 + \xi^2\right)\\
&= \xi \pi \operatorname{cosech}(\xi \pi) \prod_{j=1}^{m} \left(j^2 + \xi^2\right),
\end{align*}
and the formula is proven.
\end{proof}

Now apply the Poisson summation formula in order to get the desired series values.

\begin{theorem}\label{thm:even_sech}
The following identities hold:
\begin{align*}
(i) \quad \sum_{k\in\mathbb{Z}} \mathrm{sech}^2 (k \pi) =&\ G^2 + \frac{1}{\pi},\\
(ii) \quad \sum_{k \in \mathbb{Z}} \mathrm{sech}^4(k \pi) =& \frac{2}{3} \left(G^4 + G^2 + \frac{1}{\pi}\right) \\
(iii) \quad \sum_{k \in \mathbb{Z}} \mathrm{sech}^6(k \pi) =&\ \frac{2}{5} G^6 + \frac{2}{3} G^4 + \frac{8}{15}\left(G^2 + \frac{1}{\pi}\right) \\
(iv) \quad \sum_{k \in \mathbb{Z}} \mathrm{sech}^8(k \pi) =&\ \frac{8}{35} G^8 + \frac{8}{15} G^6 + \frac{28}{45} G^4 + \frac{16}{35}\left(G^2 + \frac{1}{\pi}\right)
\end{align*}
with the Gau\ss ian constant $G = \varpi/\pi$.
\end{theorem}

\begin{proof}
\textit{Ad (i)}: For the application of the Poisson summation formula, we first calculate the Fourier transformation of $\operatorname{sech}^2(k\pi)$. According to Lemma \ref{le:fourier_sech}, one has
$$
\mathcal{F}\left[x \mapsto \operatorname{sech}^2(x \pi)\right](\xi) =  \frac{2}{\pi} \left|\Gamma(1 + i\xi) \right|^2 = 2\xi \operatorname{cosech}(\xi \pi).
$$
Now, the Poisson summation formula and Equation \eqref{eq:sum_cosech} in the proof of Proposition \ref{prop:lemniscatic_sech} yield
$$
\sum_{\kappa \in \mathbb{Z}} \operatorname{sech}^2(\kappa \pi) = \sum_{\kappa\in\mathbb{Z}} \mathcal{F}\left[x \mapsto \operatorname{sech}^2(x \pi)\right](\xi)\big|_{\xi = \kappa} = \sum_{\kappa \in \mathbb{Z}} 2\xi \operatorname{cosech}(\xi \pi)\big|_{\xi = \kappa} = G^2 + \frac{1}{\pi},
$$
as claimed.\\

\textit{Ad (ii)} [Probabilistic proof]: Let $X\sim\mathsf{Cosech}$ be a discrete random variable which is distributed according to the hyperbolic cosecant distribution with normalization constant $C=\pi/(\pi G^2 + 1)$. Therefore, by using the duality of the Fourier transformation from \textit{(i)}, the Poisson summation formula yields
\begin{align*}
\mathbb{E}\left[X^2\right] &= \frac{\pi}{\pi G^2 + 1} \sum_{\kappa \in \mathbb{Z}} \kappa^2 \cdot 2\kappa \operatorname{cosech}(\kappa \pi)\\
&= \frac{\pi}{\pi G^2 + 1} \sum_{\kappa \in \mathbb{Z}} \mathcal{F}\left[x \mapsto x^2 \cdot 2x \operatorname{cosech}(x \pi)\right](\xi)\big|_{\xi = \kappa}\\
&= \frac{\pi}{\pi G^2 + 1} \sum_{\kappa \in \mathbb{Z}} \left(\frac{i}{2\pi}\right)^2 \frac{\mathrm{d}^2}{\mathrm{d}\xi^2} \mathcal{F}\left[x \mapsto 2x \operatorname{cosech}(x \pi)\right](\xi)\big|_{\xi = \kappa}\\
&= \frac{\pi}{\pi G^2 + 1} \sum_{\kappa \in \mathbb{Z}} \frac{-1}{4\pi^2} \frac{\mathrm{d}^2}{\mathrm{d}\xi^2} \operatorname{sech}^2(\xi \pi)\big|_{\xi = \kappa}\\
&= \frac{\pi}{\pi G^2 + 1} \left( \frac{3}{2} \sum_{\kappa \in \mathbb{Z}} \operatorname{sech}^4(\kappa \pi) - \sum_{\kappa \in \mathbb{Z}} \operatorname{sech}^2(\kappa \pi) \right).\\
\end{align*}

But Proposition \ref{prop:Var_cosech} yields $\mathbb{E}\left[X^2\right] = \operatorname{Var}(X) = \pi G^4/(\pi G^2 + 1)$ and from \textit{(i)} we get $\sum_{k\in\mathbb{Z}} \mathrm{sech}^2 (k \pi) = G^2 + 1/\pi$. Thus, the equations above hold if and only if
$$
\sum_{\kappa \in \mathbb{Z}} \operatorname{sech}^4(\kappa \pi) = \frac{2}{3}\left(G^4 + G^2 + \frac{1}{\pi}\right),
$$
as claimed.\\

\textit{Ad (ii)} [Analytical proof]: First, we compute the Fourier transformation of $\operatorname{sech}^4(k\pi)$. Lemma \ref{le:fourier_sech} yields
$$
\mathcal{F}\left[x \mapsto \operatorname{sech}^4(x \pi)\right](\xi) = \frac{4}{3\pi} \left|(1 + i \xi)\Gamma(1 + i\xi) \right|^2 = \frac{2}{3} (1 + \xi^2) \cdot 2\xi \operatorname{cosech}(\xi \pi).
$$
Now, our Equation \eqref{eq:sum_cosech} in the proof of Proposition \ref{prop:lemniscatic_sech} provides
$$
\sum_{\kappa \in \mathbb{Z}} 2\kappa \operatorname{cosech}(\kappa \pi) = G^2 + \frac{1}{\pi}.
$$
Further, Ramanujan's series identity (i) of Entry 15 in Chapter 17 of his second notebook, cf. \cite[p.\,132]{ramanujan_1991}, which we already used in the proof of Proposition \ref{prop:Var_cosech} with $x=1/2$, $y=\pi$, and $z=\sqrt{2}G$, yields
$$
\sum_{n=1}^\infty \frac{n^{3}}{\sinh(y n)} = \frac{1}{8} z^4 x = \frac{1}{2} \sum_{\kappa\in\mathbb{Z}} \kappa^3 \operatorname{cosech}(\kappa \pi) = \frac{G^4}{4}.
$$
These formulae and the Poisson summation formula finally yield
\begin{align*}
\sum_{\kappa \in \mathbb{Z}} \operatorname{sech}^4(\kappa \pi) &= \sum_{\kappa\in\mathbb{Z}} \mathcal{F}\left[x \mapsto \operatorname{sech}^4(x \pi)\right](\xi)\big|_{\xi = \kappa}\\
&= \frac{2}{3} \left(\sum_{k \in \mathbb{Z}} 2\xi \operatorname{cosech}(\xi  \pi)\big|_{\xi = \kappa} + 2 \sum_{\kappa\in\mathbb{Z}} \xi^3 \operatorname{cosech}(\xi \pi)\big|_{\xi = \kappa}\right)\\
&= \frac{2}{3} \left(G^4 + G^2 + \frac{1}{\pi}\right),
\end{align*}
as claimed.\\

\textit{Ad (iii)}: Recall that for $X \sim \mathsf{Cosech}$ with normalization constant $C = \pi/(\pi G^2 + 1)$, we have
\[
\sum_{\kappa \in \mathbb{Z}} 2\kappa \operatorname{cosech}(\pi \kappa) = \frac{1}{C} = G^2 + \frac{1}{\pi},
\]
and for the even moments ($j \in \{2, 4, 6\}$):
\[
\sum_{\kappa \in \mathbb{Z}} \kappa^j \cdot 2\kappa \operatorname{cosech}(\pi \kappa) = \frac{1}{C} \mathbb{E}\left[X^j\right].
\]

For $m=3$, the product in Lemma \ref{le:fourier_sech} evaluates to
\[
\prod_{j=1}^{2} \left(j^2 + \xi^2\right) = \left(1 + \xi^2\right)\left(4 + \xi^2\right) = 4 + 5\xi^2 + \xi^4,
\]
with the prefactor $2^4/5! = 2/15$. Thus,
\[
\mathcal{F}\left[x \mapsto \operatorname{sech}^6(\pi x)\right](\xi) = \frac{2}{15} \left(4 + 5\xi^2 + \xi^4\right) \cdot 2\xi \operatorname{cosech}(\pi \xi).
\]
Applying the Poisson summation formula and substituting $\mathbb{E}[X^2]/C = G^4$ and $\mathbb{E}[X^4]/C = 3G^6$ yield
\begin{align*}
\sum_{k \in \mathbb{Z}} \operatorname{sech}^6(k \pi)
&= \frac{2}{15} \sum_{\kappa \in \mathbb{Z}} \left(4 + 5\kappa^2 + \kappa^4\right) \cdot 2\kappa \operatorname{cosech}(\pi \kappa) \\
&= \frac{2}{15} \cdot \frac{1}{C} \left( 4 + 5\,\mathbb{E}\left[X^2\right] + \mathbb{E}\left[X^4\right] \right) \\
&= \frac{2}{15} \left( 4\left(G^2 + \frac{1}{\pi}\right) + 5 G^4 + 3 G^6 \right) \\
&= \frac{2}{5} G^6 + \frac{2}{3} G^4 + \frac{8}{15}\left(G^2 + \frac{1}{\pi}\right),
\end{align*}
as claimed.\\

\textit{Ad (iv)}: For $m=4$, the product in Lemma \ref{le:fourier_sech} expands to
\[
\prod_{j=1}^{3} \left(j^2 + \xi^2\right) = \left(1 + \xi^2\right)\left(4 + \xi^2\right)\left(9 + \xi^2\right) = 36 + 49\xi^2 + 14\xi^4 + \xi^6,
\]
with the prefactor $2^6/7! = 4/315$. Hence,
\[
\mathcal{F}\left[x \mapsto \operatorname{sech}^8(\pi x)\right](\xi) = \frac{4}{315} \left(36 + 49\xi^2 + 14\xi^4 + \xi^6\right) \cdot 2\xi \operatorname{cosech}(\pi \xi).
\]
Applying Poisson summation along with $\mathbb{E}[X^6]/C = 18G^8$ results in
\begin{align*}
\sum_{k \in \mathbb{Z}} \operatorname{sech}^8(k \pi)
&= \frac{4}{315} \cdot \frac{1}{C} \left( 36 + 49\,\mathbb{E}\left[X^2\right] + 14\,\mathbb{E}\left[X^4\right] + \mathbb{E}\left[X^6\right] \right) \\
&= \frac{4}{315} \left( 36\left(G^2 + \frac{1}{\pi}\right) + 49 G^4 + 14 \cdot 3 G^6 + 18 G^8 \right) \\
&= \frac{4}{315} \left( 18 G^8 + 42 G^6 + 49 G^4 + 36\left(G^2 + \frac{1}{\pi}\right) \right) \\
&= \frac{8}{35} G^8 + \frac{8}{15} G^6 + \frac{28}{45} G^4 + \frac{16}{35}\left(G^2 + \frac{1}{\pi}\right),
\end{align*}
which finishes the proof.
\end{proof}

\begin{remark}
Some results of the present article and of \cite{pleschberger_2026} may also be useful in investigating the conjectured identities in \cite{vigren_2020} concerning electrostatic lattice sums.
\end{remark}


\begin{thebibliography}{9}

\bibitem[Apostol, 1990]{apostol_1990}
Apostol, Tom M.:
\textit{Modular Functions and Dirichlet Series in Number Theory (2nd edition)}. Springer.

\bibitem[Borwein and Borwein, 1987]{borwein_1987}
Borwein, Jonathan M.; Borwein, Peter B.: \textit{Pi and the AGM. A Study in Analytic Number Theory and Computational Complexity}. Wiley.

\bibitem[Bruckman, 1977]{bruckman_1977}
Bruckman, P. S.: \textit{On the evaluation of certain infinite series by elliptic functions}. Fibonacci Quart. 15 (1977), 293--310.

\bibitem[Chowla and Selberg, 1967]{selberg_1967}
Chowla, S., and Selberg, A. (1967). \textit{On Epstein's Zeta-function}. Journal für die reine und angewandte Mathematik, 227, 86--110.

\bibitem[Ling, 1974]{ling_1974}
Ling, C. B.: \textit{On Summation of Series of Hyperbolic Functions}. SIAM J. Math. Anal., 5 (1974), pp. 551--562.

\bibitem[Ling, 1988]{ling_1988}
Ling, C. B. \textit{Extension of summation of series of hyperbolic functions}. J. Math. Anal. Appl. 131(2) (1988), 451--464.

\bibitem[Jacobi, 1829]{jacobi_1829}
Jacobi, C. G. J. (1829). \textit{Fundamenta nova theoriae functionum ellipticarum}.
Publisher Sumtibus fratrum Borntraeger

\bibitem[NIST, 2026]{nist_2026}
National Institute of Standards and Technology (NIST). \textit{Digital Library of Mathematical Functions (DLMF)}. Version 1.2.8; Release date: September 2026.

\bibitem[Pleschberger, 2026]{pleschberger_2026}
Pleschberger, Leonard J. E. (2026): \textit{Discrete Hyperbolic Secant Distributions}. Via: https://arxiv.org/pdf/2609.11293.

\bibitem[Ramanujan and Berndt, 1991]{ramanujan_1991}
Ramanujan, Srinivasa; Berndt, Bruce (1991).
\textit{Ramanujan's Notebooks. Part III.}
Springer.

\bibitem[Vigren and Dieckmann, 2020]{vigren_2020}
Vigren, Erik; Dieckmann, Andreas (2020). \textit{Simple Solutions of Lattice Sums for Electric Fields Due to Infinitely Many Parallel Line Charges}. Symmetry (2020), 12(6), 1040.

\bibitem[Whittaker and Watson, 1920]{whittaker_1920}
Whittaker, Edmund T.; Watson, George N. (1920).
\textit{A Course of Modern Analysis; an Introduction to the General Theory of Infinite Processes and of Analytic Functions; with an Account of the Principal Transcendental Functions (3rd Edition).}
Cambridge University Press.

\bibitem[Zucker, 1979]{zucker_1979}
Zucker, I. J.: \textit{The Summation of Series of Hyperbolic Functions}. SIAM J. Math. Anal. 10(1) (1979), 192--206.

\bibitem[Zucker, 1984]{zucker_1984}
Zucker, I. J.: \textit{Some infinite series of exponential and hyperbolic functions}. SIAM J. Math. Anal. 15 (1984), 406--413.

\end{thebibliography}
\end{document}